\documentclass[12pt]{amsart}

\usepackage{amssymb}

\usepackage{epsfig,color}

\newtheorem{thm}{Theorem}
\newtheorem{prop}{Proposition}
\newtheorem{defin}{Definition}

\newtheorem{cor}{Corollary}
\newtheorem{lem}{Lemma}
\newtheorem{rem}{Remark}

\newcommand{\g}{{\mathfrak g}}

\newcommand{\N}{{\mathbb N}}

\newcommand{\R}{{\mathbb R}}\newcommand{\RRR}{{\mathbb R}^3}\newcommand{\Rn}{{\mathbb R}^n}

\newcommand{\Z}{{\mathbb Z}}

\newcommand{\id}{\text{Id}}   
\newcommand{\iso}{\text{Iso}}   

\newcommand{\sgn}{\text{sgn}}  

\newcommand{\vol}{\text{vol}}

\newcommand{\tc}{\tilde{c}}

\newcommand{\tka}{\tilde{\kappa}}

\newcommand{\ka}{\kappa}
\newcommand{\la}{\lambda}

\newcommand{\si}{\sigma}\newcommand{\Si}{\Sigma}
\newcommand{\vp}{\varphi}

\begin{document}

\bibliographystyle{plain}

\title[Curvatures, volumes and norms of derivatives]
{The curvatures of regular curves and euclidean invariants of
their derivatives}

\author{Eugene Gutkin}

\address{Nicolaus Copernicus University (UMK), Chopina 12/18, Torun 87-100 and
Mathematics Institute of the Polish Academy of Sciences (IMPAN), Sniadeckich 8,
Warszawa 10, Poland} \email{gutkin@mat.umk.pl,gutkin@impan.pl}



\date{\today}

\begin{abstract}
The well known formulas express the curvature and the torsion of a curve in $\RRR$ in terms
of euclidean invariants of its derivatives. We obtain expressions of this
kind for all curvatures of curves in $\R^n$. It follows that a curve in $\R^n$
is determined up to an isometry by the norms of its $n$ derivatives. We extend
these observations to curves in arbitrary riemannian manifolds. 
\end{abstract}

\maketitle

\tableofcontents

\section{Introduction}       \label{intro}
What mathematical material is better known than the curvature and the torsion
of spatial curves? It is a must in every textbook on differential
geometry. See, for instance, \cite{DoC76,Sp79,Kl78}. Although most
mathematicians are familiar with this material, I will briefly
outline it now. Let $c(t),\,a\le t \le b,$ be a differentiable
curve in $\RRR$. Orthonormalizing the derivatives
$c'(t),c''(t)$,\footnote{We assume that the vector function $c(t)$
is differentiable as many times as needed, and that the
derivatives $c'(t),c''(t)$ are linearly independent. In what
follows, we refer to conditions of this kind as the regularity
assumptions.} we associate with the curve an orthonormal triple
$e(t)=(e_1(t),e_2(t),e_3(t))$. Differentiating these vectors, we
obtain a system of linear differential equations; it is customary
to write it as $e'(t)=||c'(t)||F(t)e(t)$. The $3\times 3$ matrix
$F(t)$ is determined by the curve; it has a very special form.

These observations were obtained independently and
simultaneously\footnote{Around 1850.} by  Frenet and Serret. See
the Wikipedia article \cite{Wi} for this. The standard terminology
is as follows: $e(t)=(e_1(t),e_2(t),e_3(t))$ is the {\em
Frenet-Serret frame}, the equation $e'(t)=||c'(t)||F(t)e(t)$ is
the {\em Frenet-Serret equation}, the matrix $F(t)$ is the {\em
Frenet-Serret matrix}. It is skew-symmetric and tri-diagonal. The
above-diagonal entries of $F(t)$ are the {\em curvature} and the
{\em torsion} of the curve. These two functions, say $\ka(t)$ and
$\tau(t)$, plus the speed $||c'(t)||$, determine the curve
up to an orientation preserving isometry of
$\RRR$. Let $l>0$ be the length of our curve. Replacing $t$ by the
{\em arclength parameter}, we code the curve by the functions
$\ka:[0,l]\to\R_+$ and $\tau:[0,l]\to\R$.

Two questions arise. One is to consruct, for a given pair of
functions on $[0,l]$, the essentially unique curve
$c:[0,l]\to\RRR$ whose curvature and torsion are these functions.
This is equivalent to integration of the differential equation
$e'(\cdot)=F(\cdot)e(\cdot)$. The solution is not given, in
general, by an explicit formula. The other problem is to
explicitly determine the curvature and the torsion of
a given curve $c(\cdot)$ in $\RRR$. The following identities are well
known:\footnote{Textbooks on differential geometry usually give
them as exercises; see, for instance, \cite{Kl78} and
\cite{DoC76}. However, they are not in \cite{Sp79}.}
\begin{equation}    \label{curv_eq}
 \ka(t)=\frac{||c'(t)\times c''(t)||}{||c'(t)||^3},
\end{equation}
\begin{equation}    \label{torsion_eq}
 \tau(t)=\frac{\det(c'(t),c''(t),c'''(t))}{||c'(t)\times c''(t)||^2}.
\end{equation}
In the arclength parameter these equations further simplify:
$$
\ka=||c'\times c''||,\   \tau=\det(c',c'',c''')/\ka^2.
$$
Note that these identities involve the {\em cross product}\footnote{It is also called the vector product.}
 $u\times v$ of vectors in $\RRR$. This operation
is peculiar to $\RRR$. It has to do with the canonical isomorphism of $\RRR$ and 
the Lie algebra $\g(SO(3))$ of the group
of linear isometries in $\RRR$. With this isomorphism, $u\times v$ becomes the Lie bracket in $\g(SO(3))$.

There are modifications of the Frenet-Serret approach, as well as
generalizations to curves in other spaces
\cite{Gr74,Bi75,Gr83,Iv00,SaWa03}. The most straightforward is to
extend this approach to the curves in euclidean spaces $\R^n$,
where $n$ is arbitrary. Let $c(t),\,a\le t \le b,$ be a regular
curve in $\R^n$. Orthonormalizing the vectors
$c'(t),\dots,c^{(n-1)}(t)$, we obtain a moving orthonormal frame
$e(t)=(e_1(t),\dots,e_n(t))$. Differentiating it, we obtain the
system $e'(t)=||c'(t)||F(t)e(t)$ of linear differential equations.
The matrix $F(t)$ is skew-symmetric and tri-diagonal. Its matrix
elements yield $n-1$ {\em curvature functions} $\ka_r,1\le r \le
n-1$. Together with $||c'(\cdot)||$, they determine the curve up
to an orientation preserving isometry of  $\R^n$. The definition
of the top curvature $\ka_{n-1}$ differs somewhat from those of
the other $n-2$ curvatures; the function $\ka_{n-1}$ may well be
called the torsion of a curve in $\R^n$.

This material is due to C. Jordan \cite{Jo874}; see \cite{Wi} for
more information. The modern terminology does not acknowledge
Jordan's contribution: It is customary to say the Frenet-Serret frame,
the Frenet-Serret equation, etc, no matter the dimension of ambient
space.

The original goal of this work was to obtain analogs of
equation~~\eqref{curv_eq} and equation~~\eqref{torsion_eq} for all
of the curvatures of regular curves in euclidean spaces of
arbitrary dimensions. Theorem~~\ref{all_curv_thm} gives these
generalizations. See also Corollary~~\ref{crsprd_curv_cor} and
Corollary~~\ref{arcleng_cor}. They provide remarkably simple
expressions for all curvatures in terms of the volumes of
parallelepipeds spanned by the higher derivatives of the curve.
These expressions allow us to estimate the distortion of curvatures
under affine transformations.
See Theorem~~\ref{estim_thm} and Corollary~~\ref{estim_cor}.

The $n-1$ curvatures, together with the norm of the tangent
vector, give a complete set of invariants for curves in
$n$-dimensional euclidean spaces. More precisely, they determine
the parameterized curve, up to an isometry of the ambient space.
We point out a problem with this set of invariants: The curvatures
and the norm of the derivative have very different natures.
Theorem~~\ref{deriv_norm_thm} yields more natural
invariants. The norms of the derivatives up to the $n$th order
form a complete set of invariants for the curve in question. This result is a consequence
of Theorem~~\ref{all_curv_thm}.

From the geometry viewpoint, euclidean spaces are special
examples of riemannian manifolds. In section~~\ref{rieman} we
extend the above observations to arbitrary riemannian
manifolds. These generalisations are straightforward. We obtain
them by replacing the differentiation in $\R^n$ by the riemannian 
covariant differentiation. See Theorem~~\ref{riem_all_curv_thm},
Theorem~~\ref{Riem_deriv_norm_thm}, and
Theorem~~\ref{nonoren_main_thm}.

\section{Heuristics; connections to mathematical physics}       \label{set}
Let us try to guess the $n$-dimensional versions of
equations~~\eqref{curv_eq} and~~\eqref{torsion_eq}. Let $u,v\in\RRR$.
It is immediate from the definition that
$||u\times v||$ is equal to the area of the parallelogram
$P(u,v)\subset\RRR$ spanned by the vectors $u,v$. Let $k\le n$ and
let $v_1,\dots,v_k\in\R^n$ be any vectors. Denote by
$P(v_1,\dots,v_k)\subset\R^n$ the $k$-dimensional parallelepiped
spanned by $v_1,\dots,v_k$ and by $\vol(v_1,\dots,v_k)$ its
$k$-volume.\footnote{Note that $\vol(v_1,\dots,v_k)=0$ iff the
vectors are linearly dependent.} With this notation, we rewrite
equation~~\eqref{curv_eq} and equation~~\eqref{torsion_eq},
respectively, as
$$
\ka(t)=\frac{\vol(c'(t),c''(t))}{\vol(c'(t))^3}
$$
and
$$
\tau(t)=\frac{\det(c'(t),c''(t),c'''(t))}{\vol(c'(t),c''(t))^2}.
$$

Let $c^{(k)}=c^{(k)}(t)$ denote the $k$-th derivative of the
function $c(\cdot)$. The above expressions for $\ka$ and $\tau$ do
not explicitly contain cross products. They suggest that for
$r<n-1$ the curvature $\ka_r$ of a curve $c(\cdot)\in\R^n$ should
be expressed in terms of $i$-dimensional volumes
$\vol(c',\dots,c^{(i)})$ with $i\le r$. They  also suggest that
the formula for the torsion $\ka_{n-1}$ should contain the
determinant $\det(c',\dots,c^{(n)})$ as a factor. Moreover, the
above expressions lead one to speculate that $\ka_r$ might be a
product of powers of $\vol(c',\dots,c^{(i)})$ and
$\det(c',\dots,c^{(n)})$.

The actual expressions for curvatures given by
Theorem~~\ref{all_curv_thm} do agree with these heuristics. The
author doubts, however, that any one would guess the strikingly
simple identities in Theorem~~\ref{all_curv_thm} and
Corollary~~\ref{crsprd_curv_cor} solely from
equations~~\eqref{curv_eq} and~~\eqref{torsion_eq}.

\medskip

Our proof of Theorem~~\ref{all_curv_thm} in section~~\ref{exter}
is based on the notion of multiple cross products of vectors in a
euclidean space of any dimension. This notion is not new. See, for
instance, \cite{Sp79} for the cross product of $n-1$ vectors in
$\R^n$. In particular, multiple cross products are used in
mathematical physics. Thus, the work \cite{DaGu95} explores the
triple cross product of vectors in $\R^4$ to analyze the
generalized Heisenberg ferromagnet. We will now briefly survey the
relevant material.

The classical\footnote{As opposed to quantum.} Heisenberg
model\footnote{More precisely, the classical, isotropic Heisenberg
ferromagnet.} is described by the differential equation
$S_t=S\times S_{xx}$, where $S(x,t)$ is a differentiable function
with values in $\RRR$. It is immediate that
$||S(x,t)||$ does not depend on $t$. In the physical interpretation,
$S(x,t)$ is the spin at time $t$ located at the point $x\in\R$.
Spins are unit vectors in $\RRR$, thus $||S(x,t)||=1$. Hence, for
every $t\in\R$, we have a continuous spin chain $S(\cdot,t)$; its
time evolution is described by $S_t=S\times S_{xx}$.  We view $x$
as the arclength parameter for a time-dependent curve
$c(x,t)\in\RRR$ such that $S(x,t)=c_x(x,t)$. Then the equation
$S_t=S\times S_{xx}$ defines a time evolution for curves in
$\RRR$ parameterized by arclength. This evolution can be described by certain nonlinear
partial differential equations on the curvature and the torsion of
the curve. Besides being of interest on its own, the equation
$S_t=S\times S_{xx}$ is equivalent to the classical nonlinear Schroedinger
equation \cite{DoSa94}.\footnote{I thank A. Veselov for pointing this out to me.} 
It is not known whether the quantum Heisenberg model
is  equivalent to the quantum nonlinear Schroedinger
equation \cite{Gut88}.

The generalized Heisenberg ferromagnet studied in \cite{DaGu95}
corresponds to the time evolution of a spin chain $S(x,t)$  with
values in $\R^4$. It is given by the equation $S_t=S\times
S_x\times S_{xx}$, where $u\times v\times w$ is the triple cross
product in $\R^4$. Again, $||S(x,t)||$ does not depend on $t$,
and we set $||S(x,t)||=1$. Viewing $x$ as the arclength parameter
for a curve $c(x,t)$ satisfying $S=c_x$, we obtain a time
evolution for curves in $\R^4$. As in the case of the
Heisenberg model in $\RRR$, the evolution $c(x,t)$ is equivalent to a system
of nonlinear partial differential equations on the three
curvatures $\ka_1(x,t),\ka_2(x,t),\ka_3(x,t)$. See \cite{DaGu95}
for details.

We will now briefly discuss a generalization of the Heisenberg
model to spin chains with values in $\R^n$, $n\ge 4$. Consider the
equation
\begin{equation}       \label{n_heisen_eq}
S_t=S\times S'\times S''\times\cdots\times S^{(n-2)}.
\end{equation}
As before, $||S||$ does not change with the time. The
$n$-dimensional Heisenberg ferromagnet is given by
equation~~\eqref{n_heisen_eq} under the condition $||S||=1$. Set
$S(x,t)=c'(x,t)$.
Equation~~\eqref{n_heisen_eq} defines a time evolution for curves
in $\R^n$ parameterized by the arclength parameter. 
It is equivalent to a system of nonlinear partial
differential equations on the $n-1$ curvatures
$\ka_1(x,t),\cdots,\ka_{n-1}(x,t)$. The material exposed in the
body of the paper suggests an approach to invariants
of equation~~\eqref{n_heisen_eq}. The proposition below
illustrates this approach. The proof is straightforward, and we
leave it to the reader.

\begin{prop}         \label{1_curv_prop}
For $n\ge 3$ let $c(x,t)$ be a time-dependent curve in $\R^n$
satisfying equation~~\eqref{n_heisen_eq}. Let
$\ka_1(x,t),\dots,\ka_{n-1}(x,t)$ be its curvatures. {\em 1}. The curvature
$\ka_1$ does not depend on time. {\em 2}. If $n\ge 4$ then
$<c',c'''>$ does not depend on time.
\end{prop}


Concluding this section, we note  that evolutions of curves in
$\RRR$ by the curvature and torsion have applications to
turbulence and to DNA analysis \cite{RiKi02}.

\section{Multiple cross products  for euclidean spaces of arbitrary dimensions}   \label{exter}
By a euclidean space we will mean a finite dimensional, oriented, real vector space with a
positive definite scalar product $<\cdot,\cdot>$.
Let $V^n$ be such a space. Choosing a positive orthonormal basis, say $e_1,\dots,e_n$,
we identify the space with $\Rn=\{(x_1,\dots,x_n)\}$; then $<x,y>=x_1y_1+\dots+x_ny_n$.
Sometimes it will be convenient to use a positive orthonormal basis.
However, our approach is coordinate free. Neither our results nor our methods depend on a particular basis.

The {\em exterior algebra}
$$
\bigwedge V =\oplus_{k=0}^n\wedge^kV
$$
is endowed with several structures. First of all, each $\wedge^kV$ is a real vector space
and $\dim_{\R}\wedge^kV= {n\choose k}$. The subspaces $\wedge^kV,0\le k \le n,$
provide a grading of $\bigwedge V$; we refer to $w\in\wedge^kV$ as elements of degree $k$.
The wedge product is anticommutative: Let $z,w\in\bigwedge V$ have degrees $k,l$ respectively;
then
$$
w\wedge z=(-1)^{kl} z\wedge w.
$$
The pairing $<\cdot,\cdot>$ on $V$ induces a  bilinear form on
$\bigwedge V$. We denote it by $<\cdot,\cdot>$ as well. The subspaces $\wedge^kV$ are pairwise orthogonal with respect
to $<\cdot,\cdot>$; thus, we only need to determine
$<\cdot,\cdot>$ on each $\wedge^kV$. The vector space $\wedge^kV$ is spanned by elements $v_1\wedge\cdots\wedge v_k$;
thus, it suffices to define the scalar product for monomials. Let $S_k$ be the permutation group of $k$ items.
We code permutations $g\in S_k$ by $k$-tuples $(i_1,\dots,i_k)$ of distinct elements in $\{1,\dots,k\}$.
Let $\si(g)\in\{0,1\}$ be the {\em parity} of permutation, so that $\det g=(-1)^{\si(g)}$.
Then we have
\begin{equation}    \label{scal_prod_eq}
<u_1\wedge\cdots\wedge u_k,v_1\wedge\cdots\wedge v_k>=\sum_{g\in S_k}(-1)^{\si(g)}<u_1,v_{i_1}>\cdots<u_k,v_{i_k}>.
\end{equation}

Let $I=\{1\le i_1<\cdots<i_k\le n\}$
be a subset in $\{1,\dots,n\}$,  with $|I|=k$. Set $e_I=e_{i_1}\wedge\cdots\wedge e_{i_k}\in\wedge^kV$.
As $I$ runs through the subsets of $\{1,\dots,n\}$, the vectors $e_I$ form a basis of $\bigwedge V$.
By equation~~\eqref{scal_prod_eq}, the basis $\{e_I:\,I\subset\{1,\dots,n\}\}$ is orthonormal.
 Thus, the bilinear form~~\eqref{scal_prod_eq}
yields a scalar product in $\bigwedge V$.\footnote{Note that
$\bigwedge V$ does not have a natural orientation.} We will now explain why this is a natural
scalar product.

Let $v_1,\dots,v_k\in V$ be any $k$-tuple. The wedge product $v_1\wedge\cdots\wedge v_k\in\wedge^kV$
corresponds to the $k$-dimensional parallelepiped $P(v_1,\dots,v_k)$ spanned by $v_1,\dots,v_k$. Note
that the vectors $v_1,\dots,v_k$ are linearly dependent iff $P(v_1,\dots,v_k)$ collapses.
Let $\vol(v_1,\dots,v_k)$ be the $k$-volume of $P(v_1,\dots,v_k)$. We leave it to the reader to
prove the identity\footnote{For instance, by induction on $k$.}
\begin{equation}    \label{prlpd_vol_eq}
<v_1\wedge\cdots\wedge v_k,v_1\wedge\cdots\wedge v_k>=\vol(v_1,\dots,v_k)^2.
\end{equation}
Thus, the scalar product equation~~\eqref{scal_prod_eq} is the symmetric
bilinear form corresponding to the quadratic form  $\vol(v_1,\dots,v_k)^2$.
For $k=2$, equation~~\eqref{prlpd_vol_eq} yields the classical formula for the area
of a parallelogram.

Note that we have not yet used the orientation of $V$. Since $\dim\wedge^nV=1$,
the euclidean space $\wedge^nV$ is isomorphic to $\R$. There are exactly {\em two
linear isometries} $O:\R\to\wedge^nV$. Choosing one of them is equivalent to endowing $V$
with an orientation. Indeed, the space $\wedge^nV$ has two elements of unit norm.
Let $o\in\wedge^nV$ be one of them, and set $O(1)=o$. Let now
$e_1,\dots,e_n$ be an orthonormal basis in $V$. By equation~~\eqref{prlpd_vol_eq},
$||e_1\wedge\cdots\wedge e_n||=1$, hence $e_1\wedge\cdots\wedge e_n=\pm o$.
The basis $e_1,\dots,e_n$ is positive if $e_1\wedge\cdots\wedge e_n=o$, and negative otherwise.

Let now $v_1,\dots,v_n\in V$ be any vectors. Then $v_1\wedge\cdots\wedge v_n=f(v_1,\dots,v_n)o$
where $f$ is a $n$-linear form on $V$. Let $A=A(v_1,\dots,v_n)$ be the $n\times n$ matrix
of coefficients of $v_1,\dots,v_n$ with respect to any positive orthonormal basis. Then
$f(v_1,\dots,v_n)=\det A$. Thus, $\det A$ does not depend on the choice of a positive orthonormal basis.
It depends on the orientation of $V$. We set $\det A(v_1,\dots,v_n)=\det(v_1,\dots,v_n)$. Then
\begin{equation}    \label{determi_eq}
v_1\wedge\cdots\wedge v_n=\det(v_1,\dots,v_n)o.
\end{equation}

For  $z\in\bigwedge V$ let $E_z:\bigwedge V\to\bigwedge V$ be the operator of left exterior
multiplication, i. e., $E_zw=z\wedge w$. The operator of left {\em interior
multiplication} $I_z:\bigwedge V\to\bigwedge V$ is the adjoint of $E_z$ with respect to the scalar
product $<\cdot,\cdot>$, i. e., $I_z=E_z^*$. For any $z,t,w\in\bigwedge V$ we have
\begin{equation}    \label{int_mult_eq}
<I_zt,w>=<E_z^*t,w>=<t,E_zw>=<t,z\wedge w>.
\end{equation}
\begin{defin}      \label{poincare_def}
{\rm For $z\in\bigwedge V$ set $D(z)=I_zo$. We call the linear operator
$D:\bigwedge V\to\bigwedge V$ the {\em poincare duality operator}.
}
\end{defin}

It will be sometimes convenient to write
$E(z),I(w)$ and $Dt$ or $D\cdot t$ for $E_z,I_w$ and $D(t)$ respectively. The following lemma
summarizes the basic properties of these  operators.

\begin{lem}      \label{poincare_lem}
{\em 1}. Let $D':\bigwedge V\to\bigwedge V$ be the poincare duality operator
corresponding to $V$ with the orientation reversed. Then $D'=-D$.

\noindent{\em 2}. Let $u,v\in V$. Then
\begin{equation}      \label{relate_eq}
E_uI_v+I_vE_u=<u,v>\id.
\end{equation}

\noindent{\em 3}. The operator $D$ is an isometry of $\bigwedge V$.

\noindent{\em 4}. We have
$$
D^2|_{\wedge^kV}=(-1)^{k(n-k)}\id.
$$
\begin{proof}
Reversing the orientation of $V$ is equivalent to replacing the element $o\in\wedge^nV$ by $-o$.
Thus, claim 1 follows from Definition~~\ref{poincare_def}.

Let $u_1,\dots,u_k\in V$ and $v_1,\dots,v_k\in V$ be arbitrary sequences of $k$ vectors. We denote
by $G=G(u_1,\dots,u_k;v_1,\dots,v_k)$ the $k\times k$ matrix such that $G_{i,j}=<u_i,v_j>$.
Thus, $G$ is the Gram matrix corresponding to $u_1,\dots,u_k$ and $v_1,\dots,v_k$.
It is immediate from equation~~\eqref{scal_prod_eq} that
\begin{equation}      \label{det_eq}
 \det G = <u_1\wedge\cdots\wedge u_k,v_1\wedge\cdots\wedge v_k>.
\end{equation}

We will use the following notational conventions. By $v_1\wedge\cdots\wedge\widehat{v_i}\wedge\cdots\wedge v_k$
we indicate that the factor $v_i$ is omitted. By $v_1\wedge\cdots\wedge(v_i\mapsto u)\wedge\cdots\wedge v_k$
we indicate that the factor $v_i$ is replaced by $u$. Let $u,v,z_1,\dots,z_k\in V$ be arbitrary.
From the definition of operators $I_w$ and equation~~\eqref{det_eq}, we obtain
\begin{equation}      \label{ext_prod_eq}
I_u(z_1\wedge\cdots\wedge z_k) = \sum_{i=1}^k(-1)^{i-1}<u,z_i>z_1\wedge\cdots\wedge\widehat{z_i}\wedge\cdots\wedge z_k.
\end{equation}
From equation~~\eqref{ext_prod_eq} we straightforwardly calculate
$$
E_uI_v(z_1\wedge\cdots\wedge z_k) = \sum_{i=1}^k<v,z_i>z_1\wedge\cdots\wedge(z_i\mapsto u)\wedge\cdots\wedge z_k
$$
and
$$
I_vE_u(z_1\wedge\cdots\wedge z_k) =
$$
$$
<u,v>z_1\wedge\cdots\wedge z_k-\sum_{i=1}^k<v,z_i>z_1\wedge\cdots\wedge(z_i\mapsto u)\wedge\cdots\wedge z_k.
$$
Claim 2 follows.

We will now prove that $D:\bigwedge V\to\bigwedge V$ is an isometry. Since $\bigwedge V=\oplus_{k=0}^n\wedge^kV$,
it suffices to show that $D|_{\wedge^kV}$ is an isometry for $0\le k \le n$. For $k=0$ this is immediate from the definition.
Let $u,v\in V$. By claim 2, we have
$$
<Du,Dv>=<I_uo,I_vo>=<o,E_uI_vo>=
$$
$$
<u,v><o,o>-<o,I_vE_uo>=<u,v>.
$$
Thus, $D|_{V}$ is an isometry. Let now $u_1,\dots,u_k$ and $v_1,\dots,v_k$ be arbitrary vectors in $V$.
Iterating the above procedure, and using equation~~\eqref{relate_eq} every time we switch the order of
operators $E_u,I_v$, we prove by induction on $k$ that
$$
<D(u_1\wedge\cdots\wedge u_k),D(v_1\wedge\cdots\wedge v_k)>=
$$
$$
\sum_{g\in S_k}(-1)^{\si(g)}<u_1,v_{i_1}>\cdots<u_k,v_{i_k}>.
$$
Since $\wedge^kV$ is spanned by monomials, and in view of equation~~\eqref{scal_prod_eq},
this proves claim 3.

Let
$$
(z,w)=<Dz,w>=<o,z\wedge w>
$$
be the bilinear form on $\bigwedge V$ corresponding to the operator $D$. The subspaces $\wedge^iV,\wedge^jV$
are orthogonal with respect to $(\cdot,\cdot)$ unless $i+j=n$. By equation~~\eqref{determi_eq}, the bilinear
form pairs up $\wedge^kV$ and $\wedge^{n-k}V$ for all $k$. Let $z\in\wedge^kV,w\in\wedge^{n-k}V$. It is immediate from
the definition of $(\cdot,\cdot)$ and the anticommutativity of the wedge product that
$$
(z,w)=(-1)^{k(n-k)}(w,z).
$$
Denote by $D^*$ the adjoint operator with respect to $<,>$. Then for any $z\in\wedge^kV,w\in\wedge^{n-k}V$
we have
$$
(z,w)=<Dz,w>=<z,D^*w>=(-1)^{k(n-k)}(w,z)=
$$
$$
(-1)^{k(n-k)}<Dw,z>=(-1)^{k(n-k)}<z,Dw>.
$$
Thus, $D^*|_{\wedge^{n-k}V}=(-1)^{k(n-k)}D|_{\wedge^{n-k}V}$. On the other hand, by claim 3, $D^*=D^{-1}$. Hence,
$D^{-1}|_{\wedge^{n-k}V}=(-1)^{k(n-k)}D|_{\wedge^{n-k}V}$, which proves claim 4.
\end{proof}
\end{lem}
\begin{defin}      \label{cross_prod_def}
Let $v_1,\dots,v_k\in V$ be arbitrary vectors, and let $1\le k \le n$. We define the
{\em cross product} on $k$ factors
$v_1\times\cdots\times v_k$ by
\begin{equation}    \label{cross_prod_eq}
v_1\times\cdots\times v_k=D(v_1\wedge\cdots\wedge v_k).
\end{equation}
\end{defin}

Thus, the cross product on $k$ factors is a $k$-linear map from $V$ to $\wedge^{n-k}V$.
Let $k=n-1$. Then $v_1\times\cdots\times v_{n-1}\in V$. In the special case $n=3,k=2$
equation~~\eqref{cross_prod_eq} yields the standard cross product of vectors in $\RRR$.
This motivates the terminology. By definition, the operation $v_1\times\cdots\times v_k$
is distributive. However, it is not associative. In fact, already  the standard cross product
in   $\RRR$ is not associative. The following proposition summarizes the properties of cross product
that we will use.

\begin{prop}            \label{basic_cross_prop}
{\em 1}. The cross product $v_1\times\cdots\times v_k$ changes sign if we switch around any two  consecutive factors.

\noindent {\em 2}. We have
\begin{equation}       \label{crsprd_vol_eq}
||v_1\times\cdots\times v_k||=\vol(v_1,\dots,v_k).
\end{equation}

\noindent {\em 3}. Let $v_1,\dots,v_{n-1}\in V$ be linearly independent. Then the vector
$v_1\times\cdots\times v_{n-1}\in V$ has the following properties:

\noindent i) Its norm satisfies $||v_1\times\cdots\times v_{n-1}||=\vol(v_1,\dots,v_{n-1})$;

\noindent ii) The vector $v_1\times\cdots\times v_{n-1}$ is orthogonal to $v_1,\dots,v_{n-1}$;

\noindent iii) The basis $v_1,\dots,v_{n-1},v_1\times\cdots\times v_{n-1}$ is positive.

\noindent {\em 4}. Let $u_1,\dots,u_k$ and $v_1,\dots,v_{n-k}$ be arbitrary vectors in $V$. Then
\begin{equation}    \label{determ_eq}
<u_1\times\cdots\times u_k,v_1\wedge\cdots\wedge v_{n-k}>=\det(u_1,\dots,u_k,v_1,\dots,v_{n-k}).
\end{equation}
\begin{proof}
Claim 1 is immediate from the anticommutativity of the wedge product and equation~~\eqref{cross_prod_eq}.
Claim 2 follows from equation~~\eqref{cross_prod_eq}, equation~~\eqref{prlpd_vol_eq}, and
claim 3 in Lemma~~\ref{poincare_lem}. We will now prove claim 4.
From preceding equations, we have
$$
<u_1\times\cdots\times u_k,v_1\wedge\cdots\wedge v_{n-k}>=<E^*(u_1\wedge\cdots\wedge u_k)o,v_1\wedge\cdots\wedge v_{n-k}>
$$
$$
=<o, u_1\wedge\cdots\wedge u_k\wedge v_1\wedge\cdots\wedge v_{n-k}>.
$$
Equation~~\eqref{determ_eq} now follows from equation~~\eqref{determi_eq}.
In the special case $k=n-1$ equation~~\eqref{determ_eq} yields
\begin{equation}    \label{def_crss_prd_eq}
<v_1\times\cdots\times v_{n-1},v>=\det(v_1,\dots,v_{n-1},v).
\end{equation}
Claim 3 follows from equations~~\eqref{def_crss_prd_eq} and~~\eqref{crsprd_vol_eq}.
\end{proof}
\end{prop}

\begin{rem}               \label{poincare_rem}
{\em
i) We point out that equation~~\eqref{determ_eq} is equivalent to our definition of the cross product.
ii) The name we use for the operator $D$ in Definition~~\ref{poincare_def} is motivated by the following
observation. Let $T^n=\R^n/\Z^n$ be the standard torus. Set $V=H^1(T^n,\R)$. The scalar product on
$V$ is induced by the isomorphism $H_1(T^n,\R)=H^1(T^n,\R)$ and the integration with
respect to the  riemannian volume form. The orientation of $V$ comes from the orientation of $T^n$.
Then under the isomorphism $\bigwedge V = H^*(T^n,\R)$ the operator $D$ in Definition~~\ref{poincare_def}
goes to the poincare duality operator on $H^*(T^n,\R)$.

}
\end{rem}
\section{The Gram-Schmidt orthogonalization}   \label{gram_schmi}
Let $f_1,\dots,f_k\in V$ be arbitrary linear independent vectors. Then there exists a
unique collection of orthonormal vectors $e_1,\dots,e_k\in V$  such that for $1\le i \le k$ we have
\begin{equation}    \label{span_eq}
e_i=a_{i,1}f_1+\cdots+a_{i,i}f_i,\ a_{i,i}>0.
\end{equation}
The collection $e_1,\dots,e_k$ is called the Gram-Schmidt
orthogonalization of $f_1,\dots,f_k$. The coefficients in equation~~\eqref{span_eq}
are determined by the scalar products of vectors $f_1,\dots,f_k$.
We will need only the coefficients $a_{i,i}$.

\begin{lem}  \label{gram_schmi_lem}
Let $f_1,\dots,f_k\in V$ be linear independent. Let $e_1,\dots,e_k$
be the Gram-Schmidt orthogonalization. Let $a_{i,j},j\le i\le k,$ be the coefficients in equation~~\eqref{span_eq}.
Then for $1\le i \le k$ we have
\begin{equation}    \label{coeff_eq}
a_{i,i}=\frac{\mbox{vol}(f_1,\dots,f_{i-1})}{\vol(f_1,\dots,f_{i-1},f_i)}.
\end{equation}
\begin{proof}
Observe that for $i=1$ the numerator in equation~~\eqref{coeff_eq} is not defined.
By convention, the volume of the parallelepiped formed by an empty collection of vectors is $1$.
With this convention, equation~~\eqref{coeff_eq} obviously holds for $i=1$.

By the definition of Gram-Schmidt orthogonalization, for any $1 \le i \le k$ we have
$$
\vol(f_1,\dots,f_{i-1},e_i)=\vol(f_1,\dots,f_{i-1}).
$$
By equation~~\eqref{prlpd_vol_eq} and equation~~\eqref{coeff_eq}
$$
\vol(f_1,\dots,f_{i-1},e_i)=\vol(f_1,\dots,f_{i-1},f_i)a_{i,i}.
$$
Combinining the two equations, we obtain the claim.
\end{proof}
\end{lem}
\section{Universal identities for the  curvatures}   \label{all_curv}
Let $V^n$ be a euclidean space. By a {\em regular curve} in $V$ we
will mean a mapping $c:[a,b]\to V$ such that i) the interval
$[a,b]$ is nontrivial; ii) the vector function $c(t)$ is as smooth
as necessary; iii) the vectors $c'(t),\dots,c^{(n-1)}(t)$ are
linearly independent for any $t\in[a,b]$. Although it is customary
to think of the variable $t\in[a,b]$ as the time, we will denote
the differentiation with respect to $t$ by ``prime'', as opposed to
``dot''.

We will use the following notational convention. Let $E(n)$ be an expression that depends
explicitly on $n\in\N$. If the expression is defined only for $n\ge n_0$, we set
$E(k)=1$ for $k<n_0$. For instance, if $E(n)=\vol(c',c''\dots,c^{(n-1)})$, then $E(1)=1$.
We will now state and prove the main result.

\begin{thm}        \label{all_curv_thm}
Let $c(t)$ be a regular curve in $V$. Then the following holds.

\noindent{\em 1}. For $1\le r \le n-2$ the curvatures $\ka_r$ satisfy
\begin{equation}    \label{gen_curv_eq}
\ka_r=\frac{\vol(c',\dots,c^{(r-1)})\vol(c',\dots,c^{(r+1)})}{\vol(c',\dots,c^{(r)})^2}||c'||^{-1}.
\end{equation}
\noindent{\em 2}. For the torsion, i. e., the top curvature, we have
\begin{equation}    \label{top_curv_eq}
\ka_{n-1}=\frac{\vol(c',\dots,c^{(n-2)})\det(c',\dots,c^{(n)})}{\vol(c',\dots,c^{(n-1)})^2}||c'||^{-1}.
\end{equation}
\begin{proof}
Let $e(t)=(e_1(t),\dots,e_n(t))$ be the associated Frenet-Serret frame.
The Frenet-Serret equation says
\begin{equation}    \label{frenet_eq}
e'\ = \ ||c'||
\begin{pmatrix}
0        &     \ka_1    &     0         &    \cdots       &     0       \\
-\ka_1   &       0      &    \ka_2      &     \cdots      &     0       \\
\vdots   &    \ddots    &    \ddots     &     \ddots      &    \vdots        \\
0        &     \cdots   &   -\ka_{n-2}  &       0         &    \ka_{n-1}  \\
0        &     \cdots   &        0      &    -\ka_{n-1}   &      0
\end{pmatrix}
e.
\end{equation}

Let $v_1(t),\dots,v_k(t)$ be arbitrary differentiable functions with values in $V$. Then,
from Definition~~\ref{cross_prod_def}, $\left(v_1\times\cdots\times v_k\right)'
=\sum_{i=1}^kv_1\times\cdots v_i'\times\cdots v_k$.

We will refer to this identity as the {\em product rule}.
Let $1\le r <n$. By the product rule and equation~~\eqref{frenet_eq}
$$
(e_1\times\cdots\times e_r)'=||c'||\ka_r e_1\times\cdots\times e_{r-1}\times  e_{r+1}.
$$
By Lemma~~\ref{gram_schmi_lem}
\begin{equation}    \label{gram_schmi_eq}
e_1\times\cdots\times e_r=\vol(c',\dots,c^{(r)})^{-1}c'\times\cdots\times c^{(r)}.
\end{equation}
Hence, by the product rule and Proposition~~\ref{basic_cross_prop}
$$
(e_1\times\cdots\times e_r)'=
$$
$$
\left[\vol(c',\dots,c^{(r)})^{-1}\right]'c'\times\cdots\times c^{(r)}
+\vol(c',\dots,c^{(r)})^{-1}c'\times\cdots\times c^{(r-1)}\times c^{(r+1)}.
$$
We assume first that $r<n-1$ and take the scalar product with $e_r\wedge e_{r+2}\wedge\cdots\wedge e_n$.
From the former of the above equations and the product rule, we have
$$
<(e_1\times\cdots\times e_r)',e_r\wedge e_{r+2}\wedge\cdots\wedge e_n>=
$$
$$
-||c'||\ka_r\det(e_1,\dots,e_n)=-||c'||\ka_r.
$$
The latter of the above equations implies
$$
<(e_1\times\cdots\times e_r)',e_r\wedge e_{r+2}\wedge\cdots\wedge e_n>=
$$
$$
\left[\vol(c',\dots,c^{(r)})^{-1}\right]'<c'\times\cdots\times c^{(r)},e_r\wedge e_{r+2}\wedge\cdots\wedge e_n>+
$$
$$
\vol(c',\dots,c^{(r)})^{-1}<c'\times\cdots\times c^{(r-1)}\times c^{(r+1)},e_r\wedge e_{r+2}\wedge\cdots\wedge e_n>.
$$
By equation~~\eqref{gram_schmi_eq} and  equation~~\eqref{determ_eq}, the former of the two scalar products vanishes.
Using Lemma~~\ref{gram_schmi_lem} again, we obtain
$$
<(e_1\times\cdots\times e_r)',e_r\wedge e_{r+2}\wedge\cdots\wedge e_n>=
$$
$$
-\vol(c',\dots,c^{(r)})^{-1}O^{-1}(c'\wedge\dots\wedge c^{(r-1)}\wedge e_r\wedge c^{(r+1)}\wedge e_{r+2}\wedge\cdots\wedge e_n).
$$
Applying Lemma~~\ref{gram_schmi_lem} to $c^{(k+1)}$ in the above wedge product, we have
$$
<(e_1\times\cdots\times e_r)',e_r\wedge e_{r+2}\wedge\cdots\wedge e_n>=
$$
$$
-\vol(c',\dots,c^{(r)})^{-2}\vol(c',\dots,c^{(r+1)})O^{-1}(c'\wedge\dots\wedge c^{(r-1)}\wedge e_r\cdots\wedge e_n).
$$
Applying Lemma~~\ref{gram_schmi_lem} once more yields
$$
<(e_1\times\cdots\times e_r)',e_r\wedge e_{r+2}\wedge\cdots\wedge e_n>=
$$
$$
-\frac{\vol(c',\dots,c^{(r-1)})\vol(c',\dots,c^{(r+1)})}{\vol(c',\dots,c^{(r)})^2}O^{-1}(e_1\wedge\dots\wedge e_n)
$$
$$
=-\frac{\vol(c',\dots,c^{(r-1)})\vol(c',\dots,c^{(r+1)})}{\vol(c',\dots,c^{(r)})^2}.
$$
Comparing this with our previous expression for $<(e_1\times\cdots\times e_r)',e_r\wedge e_{r+2}\wedge\cdots\wedge e_n>$,
we obtain equation~~\eqref{gen_curv_eq}.

Let now $r=n-1$. As before, we compare two expressions for $<(e_1\times\cdots\times e_{n-1})',e_{n-1}>$.
Recall that, as opposed to $e_1,\dots,e_{n-1}$, the vector $e_n$ does not necessarily satisfy equation~~\eqref{span_eq}.
Instead, the vector $e_n$ is chosen so that $e_1,\dots,e_n$ form a positive orthonormal basis. Let
$$
e_n=a_{n,1}c'+\cdots+a_{n,n}c^{(n)}.
$$
The argument of Lemma~~\ref{gram_schmi_lem} allows us to calculate $a_{n,n}$; it yields
\begin{equation}    \label{last_coeff_eq}
a_{n,n}=\frac{\mbox{vol}(c'\dots,c^{(n-1)})}{\det(c'\dots,c^{(n)})}.
\end{equation}
The preceding argument for $r=n-1$ and equation~~\eqref{last_coeff_eq} yield equation~~\eqref{top_curv_eq}.
\end{proof}
\end{thm}

\section{Applications to curves in euclidean spaces}   \label{corol_appli}
We begin by exposing a few immediate consequences of Theorem~~\ref{all_curv_thm}.
\subsection{Immediate corollaries}                 \label{corol_sub}
\hfill\break Theorem~~\ref{all_curv_thm} was motivated by
equations~~\eqref{curv_eq} and~~\eqref{torsion_eq} for the
curvature and torsion of curves in $\RRR$. Our first application
of Theorem~~\ref{all_curv_thm} is to the curvatures $\ka_1,\ka_2$
for curves in arbitrary euclidean spaces.

\begin{cor}   \label{low_curv_cor}
Let $c(t)$ be a regular curve in $V^n$. If $n\ge 3$ then we have
\begin{equation}    \label{1_curv_eq}
\ka_1=\frac{\vol(c',c'')}{||c'||^3}.
\end{equation}
If $n\ge 4$, then
\begin{equation}    \label{2_curv_eq}
\ka_2=\frac{\vol(c',c'',c''')}{\vol(c',c'')^2}.
\end{equation}
\begin{proof}
When $n\ge 3$, the assumptions of claim 1 in Theorem~~\ref{all_curv_thm} hold for $\ka_1$.
Equation~~\eqref{1_curv_eq} is a special case of equation~~\eqref{gen_curv_eq}.
When $n\ge 4$, the assumptions of claim 1 in Theorem~~\ref{all_curv_thm} hold for $\ka_2$.
From equation~~\eqref{gen_curv_eq}, we have
$$
\ka_2=\frac{\vol(c')\vol(c',c'',c''')}{\vol(c',c'')^2||c'||}.
$$
Since $\vol(c')=||c'||$, we obtain equation~~\eqref{2_curv_eq}.
\end{proof}
\end{cor}

\begin{rem}   \label{n=3_rem}
Note that equation~~\eqref{curv_eq}
and equation~~\eqref{torsion_eq} are the special cases of equation~~\eqref{1_curv_eq}
and equation~~\eqref{top_curv_eq} respectively.
\end{rem}

\begin{cor}                       \label{crsprd_curv_cor}
Let $c(t)$ be a regular curve in a euclidean space $V$ of $n$ dimensions. Let $\ka_1(t),\dots,\ka_{n-1}(t)$ be its curvatures.
Then  for $1\le r \le n-2$ we have
$$
\ka_{r}=\frac{||c'\times\cdots\times c^{(r-1)}||\,||c'\times\cdots\times c^{(r+1)}||}{||c'\times\cdots\times c^{(r)}||^2}||c'||^{-1}.
$$
For the torsion we have
$$
\ka_{n-1}=\frac{||c'\times\cdots\times c^{(n-2)}||}{||c'\times\cdots\times c^{(n-1)}||^2}
\det\left(c',c'',\dots,c^{(n)}\right)||c'||^{-1}.
$$
\begin{proof}
Immediate from Theorem~~\ref{all_curv_thm} and Proposition~~\ref{basic_cross_prop}.
\end{proof}
\end{cor}

\begin{cor}   \label{arcleng_cor}
Let $c(s)$ be a regular curve in a euclidean space $V$ of $n$ dimensions parameterized by
arclength. Let $\ka_1(s),\dots,\ka_{n-1}(s)$ be its curvatures. Then for $1\le r \le n-2$
\begin{equation}    \label{arclen_curv_eq}
\ka_{r}=\frac{\vol(c',\dots,c^{(r-1)})\vol(c',\dots,c^{(r+1)})}{\vol(c',\dots,c^{(r)})^2}.
\end{equation}
The top curvature satisfies
\begin{equation}     \label{arclen_tors_eq}
\ka_{n-1}=\frac{\vol(c',\dots,c^{(n-2)})}{\vol(c',\dots,c^{(n-1)})^2}
\det\left(c',c'',\dots,c^{(n)}\right).
\end{equation}
\begin{proof}
Equation~~\eqref{arclen_curv_eq} and equation~~\eqref{arclen_tors_eq} follow from
equation~~\eqref{gen_curv_eq} and equation~~\eqref{top_curv_eq} respectively, 
via $||c'||=1$.
\end{proof}
\end{cor}

\begin{cor}   \label{volumes_cor}
{\em 1}. Let $V$ be a euclidean space of $n$ dimensions; let $c(s)$ be a regular curve in $V$ parameterized by
arclength. Let $\ka_1(s),\dots,\ka_{n-1}(s)$ be its curvatures.  Then for $1\le r \le n-1$
\begin{equation}    \label{volum_eq}
\vol(c',\dots,c^{(r)})=\ka_1^{r-1}\ka_2^{r-2}\cdots\ka_{r-2}^2\ka_{r-1}
\end{equation}
and
\begin{equation}    \label{det_curv_eq}
\det(c',\dots,c^{(n)})=\ka_1^{n-1}\ka_2^{n-2}\cdots\ka_{n-2}^2\ka_{n-1}.
\end{equation}

\noindent {\em 2}. Let $V$ be as above; let $c(t)$ be a regular curve in $V$.
Let $\ka_1(t),\dots,\ka_{n-1}(t)$ be its curvatures.  Then for $1\le r \le n-1$
\begin{equation}    \label{Volum_eq}
\frac{\vol(c',\dots,c^{(r)})}{||c'||^{r(r+1)/2}}=\ka_1^{r-1}\ka_2^{r-2}\cdots\ka_{r-2}^2\ka_{r-1}
\end{equation}
and
\begin{equation}    \label{Det_curv_eq}
\frac{\det(c',\dots,c^{(n)})}{||c'||^{n(n+1)/2}}=\ka_1^{n-1}\ka_2^{n-2}\cdots\ka_{n-2}^2\ka_{n-1}.
\end{equation}
\begin{proof}
Let $x_i\ne 0,1\le i,$ be any sequence of numbers. Set, for convenience, $x_i=1$ if $i<1$.
For $i\ge 1$ set $y_i=x_i/x_{i-1}$. Then for $k\in\N$ we have
$$
x_k=y_1\cdots y_k.
$$
Setting $y_i=1$ if $i<1$ and defining $z_j=y_j/y_{j-1}$, we have for $l\in\N$
$$
y_l=z_1\cdots z_l.
$$
The above equations yield
\begin{equation}    \label{Teleskop_eq}
x_k=z_1^kx_2^{k-1}\cdots z_k.
\end{equation}

Let $c(s)$ be as in claim 1. Set $x_i=\vol(c',\dots,c^{(i)})$. Let the sequences $y_i$ and $z_i$ be as above.
Then, by equation~~\eqref{gen_curv_eq},
$z_i=\ka_{i-1}$. Hence, Corollary~~\ref{arcleng_cor} and equation~~\eqref{Teleskop_eq}
yield equation~~\eqref{volum_eq}. Let now $c(t)$ be as in claim 2. Again, set
$x_i=\vol(c',\dots,c^{(i)})$ and define the sequences $y_i$ and $z_i$ as above.
By Theorem~~\ref{all_curv_thm}, $z_i=\ka_{i-1}||c'||$. Now equation~~\eqref{Teleskop_eq}
yields equation~~\eqref{Volum_eq}.

Equations~~\eqref{det_curv_eq} and~~\eqref{Det_curv_eq} follow the same way from
equations~~\eqref{arclen_tors_eq} and~~\eqref{top_curv_eq} respectively. We leave details to the reader.
\end{proof}
\end{cor}

\begin{rem}         \label{brauner_rem}
{\em
Equations~~\eqref{Volum_eq} and~~\eqref{Det_curv_eq} are contained in \cite{Br81}.
See problem 2 on p. 100. Since these identities are equivalent to
equations~~\eqref{gen_curv_eq} and~~\eqref{top_curv_eq} respectively,
our Theorem~~\ref{all_curv_thm} is not new. However, equations~~\eqref{gen_curv_eq} and~~\eqref{top_curv_eq}
are more direct than equations~~\eqref{Volum_eq} and~~\eqref{Det_curv_eq}; our derivation
of these identities is elementary and straightforward.\footnote{There is no information in \cite{Br81}
about the solution of problem 2.} For these reasons we feel that  Theorem~~\ref{all_curv_thm}
deserves publication.

}
\end{rem}
Now we expose some less immediate consequences of Theorem~~\ref{all_curv_thm}.
\subsection{Estimates for curvatures}                 \label{appli_sub}
\hfill\break Let $c(t)$ be a regular curve in $V^n$; let
$\ka_1(t),\dots,\ka_{n-1}(t)$ be its curvatures. Let $L:V\to V$ be
a nondegenerate linear transformation; let $\la\in V^n$. Set
$\tc(t)=Lc(t)+\la$. Let $\tka_1(t),\dots,\tka_{n-1}(t)$  be the
curvatures of $\tc$. How do they relate to
$\ka_1(t),\dots,\ka_{n-1}(t)$? If $L\in SO(V)$,\footnote{We denote
by $O(V)$ (resp. $SO(V)$) the group of (resp. orientation
preserving) linear isometries.} then $\tka_r=\ka_r$ for $1\le r\le
n-1$. Conversely, the $n-1$ curvatures determine the curve up to a
transformation $\tc(t)=Lc(t)+\la$ with $L\in SO(V)$.

Assume now that $L\notin O(V)$. Theorem~~\ref{all_curv_thm} allows
us to estimate $\tka_1(t),\dots,\tka_{n-1}(t)$. In order to state
the result, we briefly recall the notion of {\em singular values}
of matrices \cite{HoJo90}. Every $n\times n$ matrix has a
decomposition $L=U\Si V$ where $U,V\in O_n$ and $\Si$ is a
diagonal matrix with non-negative entries
$\si_1\ge\si_2\ge\cdots\ge\si_n$. Then $n$ numbers
$\si_1\ge\si_2\ge\cdots\ge\si_n\ge 0$ uniquely determined by the
matrix $L$ are its singular values.

Let $L:V^n\to V^n$ be a linear mapping. Identifying $V$ with
$\R^n$, we represent $L$ by a $n\times n$ matrix. Its singular
values do not depend on the isomorphism $V=\R^n$. Thus, we can
talk about the singular values of a linear mapping $L:V\to V$. The
subject of singular values of matrices is of use in control
theory; see, for instance, \cite{Jon97}. There are nontrivial
relationships between singular values of matrices, convex geometry
and differential geometry \cite{Gut04,GJK}.
\begin{thm}             \label{estim_thm}
Let $V^n$ be a euclidean space. Let $c(t)$ be a regular curve in
$V$. Let $\ka_1(t),\dots,\ka_{n-1}(t)$ be its curvatures. Let
$L:V\to V$ be an invertible linear mapping, let $\la\in V$ be
arbitrary and set $\tc(t)=Lc(t)+\la$. Denote by
$\tka_1(t),\dots,\tka_{n-1}(t)$  the curvatures of $\tc$.

Let $\si_1\ge\si_2\ge\cdots\ge\si_n$ be the singular values of $L$. Then for $1\le r\le n-2$ we have the bounds
\begin{equation}        \label{all_curv_bound_eq}
\frac{\si_{n-r}\si_{n-r+1}\si_{n-r+2}^2\cdots\si_n^2}{\si_1^3\si_2^2\cdots\si_r^2}\,\ka_r\le\tka_r\le
\frac{\si_1^2\cdots\si_{r-1}^2\si_r\si_{r+1}}{\si_{n-r+1}^2\cdots\si_{n-1}^2\si_n^3}\,\ka_r.
\end{equation}
For the torsion we have  the bounds
\begin{equation}        \label{top_curv_bound_eq}
\frac{\si_n^2}{\si_1^2\si_2}\,|\ka_{n-1}|
\le |\tka_{n-1}|\le \frac{\si_1^2}{\si_{n-1}\si_n^2}\,|\ka_{n-1}|.
\end{equation}
\begin{proof}
We denote by $||L||$ the {\em operator norm}, i. e.,
$$
||L||=\max_{v\in V,v\ne 0}\frac{||Lv||}{||v||}.
$$
For $0\le k \le n$ we denote by $\wedge^kL:\wedge^kV\to\wedge^kV$ the induced linear operator.
Thus, $\wedge^kL$ is the $k$th exterior power of $L$. The singular
values of $\wedge^kL$ are the $n\choose k$numbers $\si_{i_1}\cdots\si_{i_k}:1\le i_1<\cdots<i_k\le n$
listed in the decreasing order. The singular values of $L^{-1}$ are $\si_n^{-1}\ge\cdots\ge\si_1^{-1}$.
The norm of an operator is equal to its largest singular value.

Let $v_1,\dots,v_k\in V$ be arbitrary independent vectors. Then
$$
||\wedge^kL^{-1}||^{-1} \vol(v_1,\dots,v_k) \le \vol(Lv_1,\dots,Lv_k) \le ||\wedge^kL|| \vol(v_1,\dots,v_k).
$$
In view of preceding remarks, we obtain
\begin{equation}        \label{volum_bound_eq}
\si_{n-k+1}\cdots\si_n \vol(v_1,\dots,v_k) \le \vol(Lv_1,\dots,Lv_k) \le \si_1\cdots\si_k \vol(v_1,\dots,v_k).
\end{equation}
Substituting the bounds equation~~\eqref{volum_bound_eq} into the first formula of Theorem~~\ref{all_curv_thm},
we obtain equation~~\eqref{all_curv_bound_eq}. The estimate equation~~\eqref{top_curv_bound_eq} is obtained
in the same fashion from the second formula of Theorem~~\ref{all_curv_thm}. 
\end{proof}
\end{thm}

The absolute values in
equation~~\eqref{top_curv_bound_eq} are due to the fact that
$\ka_{n-1}$ is not necessarily positive.
Often we have only partial infofmation about the singular values. For instance, we may know
the norms of the matrices in question.
Using  that $\si_1=||L||,\si_n=||L^{-1}||^{-1}$, we immediately obtain from Theorem~~\ref{estim_thm}
the following statement.

\begin{cor}          \label{estim_cor}
Let the setting and the notation be as in Theorem~~\ref{estim_thm}. Then
for $1\le r\le n-2$ we have the bounds
\begin{equation}                    \label{crud_curv_bound_eq}
||L^{-1}||^{-2r} ||L||^{-2r-1}\,\ka_r \le \tka_r \le ||L||^{2r}||L^{-1}||^{2r+1}\,\ka_r.
\end{equation}
We also have
\begin{equation} \label{crud_tors_bound_eq}
||L^{-1}||^{-3}||L||^{-2}\,|\ka_{n-1}| \le  |\tka_{n-1}| \le ||L||^2 ||L^{-1}||^3\, |\ka_{n-1}|.
\end{equation}
\end{cor}

Theorem~~\ref{estim_thm} and Corollary~~\ref{estim_cor} provide
very basic estimates for the curvatures of $\tc(t)=Lc(t)$.
However, as the following remark shows, these estimates are sharp.
\begin{rem}      \label{scalar_rem}
Let $a\in\R$ be any nonzero number. Set $\tc(t)=ac(t)+\la$. Then the inequalities in
Theorem~~\ref{estim_thm} and Corollary~~\ref{estim_cor} become the identities
$$
\tka_r=|a|^{-1}\ka_r:\,1\le r \le n-2;\ |\tka_{n-1}|=|a|^{-1}|\ka_{n-1}|.
$$
\begin{proof}
For convenience of the reader, we outline a proof. Since the inequalities in Corollary~~\ref{estim_cor}
are the consequences of those in Theorem~~\ref{estim_thm}, it suffices to show that they become equalities.
We have $L=a\,\id$. We assume without loss of generality that $a>0$.
Then $||L||=a,||L^{-1}||=a^{-1}$. Equation~~\eqref{crud_curv_bound_eq} yields $a^{-1}\ka_r\le\tka_r\le a^{-1}\ka_r$.
Equation~~\eqref{crud_tors_bound_eq} becomes $a^{-1}|\ka_{n-1}|\le|\tka_{n-1}|\le a^{-1}|\ka_{n-1}|$.
\end{proof}
\end{rem}

\medskip

\subsection{Natural invariants for curves}          \label{invari_sub}
\hfill\break  Let $c(t)$ be a regular curve in $V^n$. Choosing an
orthonormal basis in $V$, we associate with the curve $n$ real
functions, $c(t)=(x_1(t),\dots,x_n(t))$; they determine the
curve. However, these functions are not intrinsically defined by
the curve; they depend on the choice of a basis in $V$. By the
Frenet-Serret equation, the $n-1$ curvatures
$\ka_1,\dots,\ka_{n-1}$ together with $||c'||$ determine the
parameterized curve $c(t)$. See equation~~\eqref{frenet_eq}. The
$n$ functions $||c'(t)||,\ka_1(t),\dots,\ka_{n-1}(t)$ are
intrinsically defined by the curve. However, this is an inhomogeneous
collection of functions. The first
member of this collection does not belong with the remaining $n-1$. 
A more homogeneous collection of functions intrinsically
defined by a curve would be
$||c'(t)||,\dots,||c^{(n)}(t)||$. Do
they determine the curve up to an isometry?

Recall that a curve $c(t)$ is regular if the vectors $c'(t),\dots,c^{(n-1)}(t)$ are linearly independent
for all $t$.

\begin{defin}      \label{str_reg_def}
A curve  $c(t)$ in $V^n$ is {\em strongly regular} if it is $n$ times continuously differentiable
and the $n$ vectors $c'(t),\dots,c^{(n-1)}(t),c^{(n)}(t)$ are linearly independent
for all $t$ in the interval of definition of the curve.
\end{defin}

We have assumed on the outset that our curves are differentiable as many times as we need. Thus, the emphasis
in Definition~~\ref{str_reg_def} is not on the existence of all $n$ derivatives but on their
linear independence. Let $v_1,\dots,v_n\in V^n$ be linearly independent. Set
$$
\sgn(v_1,\dots,v_n)=\mbox{sign}[\det(v_1,\dots,v_n)].
$$
If $v_1,\dots,v_n$ are linearly dependent, we set $\sgn(v_1,\dots,v_n)=0$.
Let $c(\cdot)$ be a regular curve in $V^n$. We define $\sgn(c(t))$ by
$$
\sgn(c(t))=\sgn(c'(t),\dots,c^{(n)}(t)).
$$
Thus, $c$ is strongly regular iff $\sgn(c(t))\equiv 1$ or $\sgn(c(t))\equiv -1$. We denote it by
$\sgn(c)$ and call it the {\em sign of the curve}.
We will say that $c(\cdot)$ is a {\em right curve} (resp. {\em left curve}) if $\sgn(c)=1$
(resp. $\sgn(c)=-1$). We will need the following lemma.

\begin{lem}        \label{det_vol_lem}
Let $v_1,\dots,v_n\in V^n$ be any vectors. Then
\begin{equation}            \label{det_vol_eq}
\det(v_1,\dots,v_n)=\sgn(v_1,\dots,v_n)\vol(v_1,\dots,v_n).
\end{equation}
\begin{proof}
If  $v_1,\dots,v_n\in V^n$ are linearly dependent, equation~~\eqref{det_vol_eq}
becomes $0=0$. Thus, we assume that $v_1,\dots,v_n\in V^n$ are linearly independent.
Setting $k=n$ in equation~~\eqref{determ_eq} and using equation~~\eqref{determi_eq}, we obtain
$$
||v_1\times\cdots\times v_n|| = |\det(v_1,\dots,v_n)|.
$$
Since $\vol(v_1,\dots,v_n)=||v_1\times\cdots\times v_n||$, the claim follows.
\end{proof}
\end{lem}

\begin{rem}       \label{det_vol_rem}
{\em

Equation~~\eqref{det_vol_eq} is well known. It is essentially equivalent to the identity
$$
\det(v_1,\dots,v_n)^2=\vol(v_1,\dots,v_n)^2
$$
which is easy to prove directly, bypassing cross products. For completeness, we outline
a proof. Choosing an orthonormal basis in $V$, we identify it with $\R^n$.
Let $A$ be the $n\times n$ matrix whose columns are the vectors $v_1,\dots,v_n$.
Then $A^tA=G(v_1,\dots,v_n)$, the Gram matrix.
Computing the determinants of these matrices, we obtain the claim.

}
\end{rem}

\begin{thm}        \label{deriv_norm_thm}
Let $c(t)$ be a strongly regular curve in $V^n$. 1. The functions $||c'(t)||,\dots,||c^{(n)}(t)||$
determine the curve up to an isometry of $V$.\footnote{We do not assume that the isometry is orientation preserving.}

\noindent 2. The functions $||c'(t)||,\dots,||c^{(n)}(t)||$ and
the number $\sgn(c)\in\{1,-1\}$ determine the curve up to an
orientation preserving isometry of $V$.
\begin{proof}
Let $I$ be a finite set of indices. Recall that
$\N=\{0,1,\dots\}$. We will say that a function, say $\psi$,  is a
linear combination of derivatives of the functions $\vp_i,i\in I,$
if $\psi=\sum_{i\in I,k\in\N}a_{i,k}\vp_i^{(k)}$, and the sum is
finite. If the right hand side in this representation is a
polynomial on variables $\vp_i^{(k)}$, we say that $\psi$ is a
differential polynomial of functions $\vp_i$.

Denote by $f_{i,j}$ the functions defined by $f_{i,j}(t)=<c^{(i)}(t),c^{(j)}(t)>$.
It will suffice to consider the indices between $1$ and $n$.
We claim that each function $f_{k,l}$ is a linear combination of derivatives of functions $f_{i,i}$
where $1\le i\le n$. By symmetry, we can assume that $k\le l$. If $k=l$, there is nothing to prove.
For $l=k+1$ the claim follows from the identity $f_{k,k}'=2f_{k,k+1}$. Let $l>k+1$. We have
$$
f_{k,l}=f_{k,l-1}'- f_{k+1,l-1}.
$$
The claim now follows by induction on $l-k$.

Let $1\le k \le n$. By equation~~\eqref{scal_prod_eq}, equation~~\eqref{prlpd_vol_eq}, and the above claim,
$\vol(c',\dots,c^{(k)})^2$ is a differential polynomial of $||c^{(i)}||^2$,
where $i=1,\dots,k$. By equation~~\eqref{gen_curv_eq}  in Theorem~~\ref{all_curv_thm},
the curvatures $\ka_1,\dots,\ka_{n-2}$ are determined by $||c'||,\dots,||c^{(n-1)}||$.
Applying Lemma~~\ref{det_vol_lem} to the vectors $c',\dots,c^{(n-1)},c^{(n)}$, we obtain the identity
\begin{equation}            \label{signum_eq}
\det(c',\dots,c^{(n)})=\sgn(c)\vol(c',\dots,c^{(n)}).
\end{equation}
By equation~~\eqref{top_curv_eq} in Theorem~~\ref{all_curv_thm} and equation~~\eqref{signum_eq},
$\ka_{n-1}$ is expressed in terms of  $||c'||,\dots,||c^{(n-1)}||,||c^{(n)}||$ and $\sgn(c)$.

By Serret-Frenet, the $n-1$ curvatures and $||c'||$ determine the
curve up to an orientation preserving isometry of $V$, yielding
claim 2. Suppose now that we know only
$||c'||,\dots,||c^{(n-1)}||,||c^{(n)}||$. The two possibilities
$\sgn(c)=\pm 1$ yield two curves
$c_+(\cdot),c_-(\cdot)$.\footnote{This is a slight abuse of
language. In fact, these are two equivalence classes of curves,
where equivalent curves differ by an orientation preserving
isometry.} Claim 1 now follows from the observation that $c_-$ and
$c_+$ differ by an orientation reversing isometry.
\end{proof}
\end{thm}

\section{Curves in arbitrary riemannian manifolds}   \label{rieman}
The approach of Frenet-Serret extends to curves in  riemannian
manifolds. We will assume that our mappings, functions, etc are
differentiable as many times as necessary. Thus, to simplify the
exposition, by a riemannian manifold $(M^n,g)$ we will mean a
connected, $C^{\infty}$ riemannian manifold of at least two
dimensions. We will now briefly recall the basic material on
riemannian geometry, referring the reader to \cite{He62} or
another textbook for details.

Let $c(t)$ be a curve in $M$ defined on an interval $I\subset\R$.
Let $c'(t)\in T_{c(t)}M$ be the tangent vector. Thus $t\mapsto
c'(t)$ maps $I$ to the tangent bundle $TM$. We assume that
$c'(t)\ne 0$. Let $D$ be the {\em operator of covariant
differentiation} with respect to $c'(t)$. Using $D$, we obtain the higher derivatives of
$c(\cdot)$. We set $c''(t)=Dc'(t),c'''(t)=Dc''(t)$, etc. We will
use the notation $c^{(k)}(t)=D^kc(t)$. The details of
Frenet-Serret approach depend on the orientability of the
manifold. We will first consider the case when $M$ is oriented.

\subsection{Curves in oriented riemannian manifolds}   \label{orient_sub}
\hfill \break Analogously to the euclidean case, we introduce the
notions of regular and strongly regular curves. Throughout this
section, $(M^n,g)$ is an oriented riemannian manifold. Whenever
this does not lead to confusion, we will simply use the notation $M$ or $M^n$.

\begin{defin}                   \label{reg_riem_def}
A parameterized curve $c(t)$ in $(M^n,g)$ is (resp. strongly)
regular if for each parameter $t\in I$ the vectors
$c'(t),\dots,c^{(n-1)}(t)\in  T_{c(t)}M$ (resp.
$c'(t),\dots,c^{(n-1)}(t),c^{(n)}(t)\in  T_{c(t)}M$) are linearly
independent.
\end{defin}

When $M$ is a euclidean space, Definition~~\ref{reg_riem_def}
reproduces Definition~~\ref{str_reg_def}. Let $c(t)$ be a geodesic
in $(M^n,g)$ parameterized by arclength.\footnote{The same
statements hold if $t$ is proportional to an arclength parameter.}
Then $c''\equiv 0$. Thus, geodesics parameterized by arclength are
not (resp. strongly) regular curves if $n>2$  (resp. $n\ge 2$).

Let $c(\cdot)$ be a regular curve in $M$. Orthonormalizing the
vectors $c'(t),\dots,c^{(n-1)}(t)\in  T_{c(t)}M$, we obtain the
Frenet-Serret frame $e(t)=(e_1(t),\dots,e_n(t))$ of the curve in
$M$. Here we have $e_i(t)\in T_{c(t)}M$. The argument pertaining
to equation~~\eqref{frenet_eq} applies verbatim and yields
\begin{equation}    \label{riem_frenet_eq}
e'(t)=D\cdot e(t)\ = \ ||c'(t)||
\begin{pmatrix}
0        &     \ka_1    &     0         &    \cdots       &     0       \\
-\ka_1   &       0      &    \ka_2      &     \cdots      &     0       \\
\vdots   &    \ddots    &    \ddots     &     \ddots      &    \vdots        \\
0        &     \cdots   &   -\ka_{n-2}  &       0         &    \ka_{n-1}  \\
0        &     \cdots   &        0      &    -\ka_{n-1}   &      0
\end{pmatrix}
e(t).
\end{equation}
In particular, just like in the euclidean case, a regular curve in
$M^n$ has $n-1$ curvatures $\ka_1(t),\dots,\ka_{n-1}(t)$; the
first $n-2$ curvatures are strictly positive. There are no
restrictions on $\ka_{n-1}$ unless $c(\cdot)$ is strongly regular.
In this case $\ka_{n-1}$ does not change sign. Equation~~\eqref{riem_frenet_eq}
allows us to obtain the counterparts
of the preceding material for curves in riemannian manifolds.

\begin{lem}                          \label{exi_uni_riem_prop}
Let $I\subset\R$ be a nontrivial interval; let $f_0,f_1,\dots,f_{n-1}$ be smooth functions
on $I$ satisfying $f_0,f_1,\dots,f_{n-2}>0$.
Let $0\in I$ be an interior point, and let $m_0\in M$ be arbitrary. Let $v_0\in T_{m_0}M$ be such that
$||v_0||=f_0(0)$. Then there exists a unique regular curve $c:I\to M$ such that

\noindent i) We have $c(0)=m_0,c'(0)=v_0$;

\noindent ii) For $t\in I$ we have $||c'(t)||=f_0(t)$ for all $t\in I$;

\noindent iii) For all $t\in I$ and $1\le i \le n-1$ we have $\ka_i(t)=f_i(t)$.
\begin{proof}
We rewrite equation~~\eqref{riem_frenet_eq} in local coordinates;
then we apply the classical propositions about the solutions of
ordinary differential equations.
\end{proof}
\end{lem}

Note that equation~~\eqref{frenet_eq} and essentially the same
argument yield the corresponding claims for curves in $\R^n$. The
only difference is that $\R^n$ has global coordinates.

We will now extend Theorem~~\ref{all_curv_thm} to curves in
riemannian manifolds. Let $c:I\to M$ be a regular curve in $M$.
Since $M$ is oriented, every tangent space $T_mM$ is a euclidean
space. Set $V(t)=T_{c(t)}M$. Then $t\mapsto V(t)$ is a smooth
function with values in $n$-dimensional euclidean spaces. Each
$V(t)$ is endowed with cross products and the other structures
defined in section~~\ref{exter}. For  $i=1,\dots,n$ let $v_i(t)\in
V(t)$ be differentiable functions. Let $1\le k \le n$. Then
$v_1(t)\times\cdots\times v_k(t)$ are differentiable vector
functions on $I$ with values in $\wedge^{n-k}V(t)$. Analogously,
$\vol(v_1(t),\dots,v_k(t))$ and $\det(v_1(t),\dots,v_n(t))$ are
differentiable real valued functions. We call them the $k$-volume,
$1\le k \le n$, and the determinant. The material in section
~~\ref{exter} straightforwardly extends to the present setting.

\begin{thm}            \label{riem_all_curv_thm}
Let $M$ be an oriented  riemannian manifold of at least two
dimensions. Let $c(t),t\in I,$ be a regular curve in $M$. For
$1\le i \le n$ let $c^{(i)}(t)\in T_{c(t)}M$ be the consecutive
covariant derivatives. Let $\vol(c',\dots,c^{(k)})$ and
$\det(c',\dots,c^{(n)})$ be the $k$-volume functions and the
determinant function.\footnote{We suppress $t$ from notation,
whenever this does not cause confusion.}  Then for $1\le r \le
n-2$ we have
\begin{equation}    \label{riem_gen_curv_eq}
\ka_r=\frac{\vol(c',\dots,c^{(r-1)})\vol(c',\dots,c^{(r+1)})}{\vol(c',\dots,c^{(r)})^2}||c'||^{-1}.
\end{equation}
The torsion, i. e., the top curvature, satisfies
\begin{equation}    \label{riem_top_curv_eq}
\ka_{n-1}=\frac{\vol(c',\dots,c^{(n-2)})\det(c',\dots,c^{(n)})}{\vol(c',\dots,c^{(n-1)})^2}||c'||^{-1}.
\end{equation}
\begin{proof}
The proof of Theorem ~~\ref{all_curv_thm} is based on the vector
calculus applied to the functions $c'(t),\dots,c^{(n)}(t)$ with
values in a euclidean space $V$. The  functions
$c'(t),\dots,c^{(n)}(t)$ have values in the variable euclidean
space $V(t)$.  All of the equations used in the proof of Theorem
~~\ref{all_curv_thm} remain valid in the present context, once we
replace the differentiation of vector functions $v(t)\mapsto
v'(t)$ by the covariant derivative $v(t)\mapsto (D\cdot v)(t)$.
Thus, the proof of Theorem ~~\ref{all_curv_thm} applies verbatim
here.
\end{proof}
\end{thm}

The propositions in section~~\ref{corol_sub} are direct
corollaries of Theorem ~~\ref{all_curv_thm}. Hence, they
straightforwardly extend to curves in oriented riemannian
manifolds. We let the reader elaborate on this remark.

\medskip

Let $c:I\to M$ be a regular curve. We define the function
$\sgn(c(t))$ on $I$ the same way we did it in
section~~\ref{invari_sub} for curves in euclidean spaces. Thus,
the only values of $\sgn(c(t))$ are $\pm 1$ and $0$. The curve
$c(\cdot)$ is strongly regular iff $\sgn(c)\equiv 1$ or
$\sgn(c)\equiv -1$. In what follows we assume without loss of
generality that $0\in I$.

\begin{prop}            \label{riem_deriv_norm_prop}
Let $M$ be an oriented riemannian manifold of at least two
dimensions. Let $m_0\in M$ and let $v_0\in T_{m_0}M$ be a nonzero
vector. Let $c:I\to M$ be a regular curve satisfying
$c(0)=m_0,c'(0)=v_0$. Then the curve $c:I\to M$ is determined by
the $n+1$ functions $||c'(t)||,\dots,||c^{(n)}(t)||$,
$\sgn(c(t))$ on $I$.
\begin{proof}
The first part of the proof of Theorem~~\ref{deriv_norm_thm} is
valid in the present setting. Hence, for $1\le k \le n$ the
$k$-volumes $\vol(c',\dots,c^{(k)})(t)$ are determined by  the
functions $||c'(t)||,\dots,||c^{(n)}(t)||$.
Equation~~\eqref{signum_eq} is also valid, thus
$\det(c',\dots,c^{(n)})(t)$ is also determined by our $n+1$
functions. The claim now follows from
Theorem~~\ref{riem_all_curv_thm} and
Lemma~~\ref{exi_uni_riem_prop}.
\end{proof}
\end{prop}

The following is immediate from
Proposition~~\ref{riem_deriv_norm_cor}.

\begin{cor}            \label{riem_deriv_norm_cor}
Let $M^n$ be an oriented riemannian manifold. Let $m_0\in M$ and
let $v_0\in T_{m_0}M$ be a nonzero vector. Let $I\subset \R$ be a
nontrivial interval containing zero. Let $f_1,\dots,f_n$ be smooth
positive functions on $I$. Let $\vp$ be a piecewise constant
function on $I$ taking values $\pm 1$.

Then there exists at most one regular curve $c:I\to M$ satisfying
i) $c(0)=m_0,c'(0)=v_0$; ii)
$||c'(t)||=f_1(t),\dots,||c^{(n)}(t)||=f_n(t)$; iii)
$\sgn(c(t))=\vp(t)$.
\end{cor}

\subsection{Curves in oriented two-point homogeneous spaces}
\label{2-point_sub}
\hfill \break In section~~\ref{corol_appli} we did not specify the
initial points and the initial directions of curves in euclidean
spaces. For instance, we stated that a curve $c(\cdot)$ in $\R^n$
whose $n-1$ curvatures and $||c'(\cdot)||$ are prescribed, is
unique up to an isometry of $\R^n$.  In order to formulate the
riemannian counterpart of the material in
section~~\ref{corol_appli}, we briefly recall the basic notions
pertaining to homogeneous riemannian manifolds.

A riemannian manifold $M$ is homogeneous if the group of
isometries $\iso(M)$ acts transitively on $M$. We denote by
$\iso_0(M)\subset\iso(M)$ the group of orientation preserving
isometries. Let $d(\cdot,\cdot)$ denote the riemannian distance.
Then $M$ is a two-point homogeneous space if for any two pairs of
points $x,y$ and $x_1,y_1$ such that $d(x,y)=d(x_1,y_1)$ there
exists $g\in\iso(M)$ satisfying $g(x)=x_1,g(y)=y_1$. Basic facts
about two-point homogeneous spaces \cite{He62} imply the
following.
\begin{lem}         \label{2_point_lem}
Let $M$ be a two-point homogeneous space. For $i=1,2$ let $m_i\in
M, v_i\in T_{m_i}M$ be such that $||v_1||=||v_2||\ne 0$. Then
there exists $g\in\iso(M)$ such that $g(m_1)=m_2,g(v_1)=v_2$.
\end{lem}

\begin{rem}         \label{2_point_rem}
{\em

The isometry claimed in Lemma~~\ref{2_point_lem} is not unique, in
general. Let $M$ be an oriented two-point homogeneous space of at
least two dimensions. Then there exist orientation preserving and
orientation reversing isometries satisfying the conditions of
Lemma~~\ref{2_point_lem}.

}
\end{rem}

\begin{prop}                          \label{uni_riem_prop}
Let $I\subset\R$ be a nontrivial interval. Let $n\ge 2$ and let
$f_0,f_1,\dots,f_{n-1}$ be smooth functions on $I$ satisfying
$f_0,f_1,\dots,f_{n-2}>0$.

Let $M^n$ be an oriented two-point homogeneous space. Then there
is a regular curve $c:I\to M$ such that for all $t\in I$ we have
i) $||c'(t)||=f_0(t)$; ii) $\ka_i(t)=f_i(t)$ for $1\le i \le n-1$.
The curve $c(t)$ is unique up to an orientation preserving
isometry of $M$.
\begin{proof}
We assume without loss of generality that $I$ contains $0$ in its
interior.  Let $m_0\in M$ be any point; let $v_0\in T_{m_0}M$ be
any vector satisfying $||v_0||=f_0(0)$. By
Lemma~~\ref{exi_uni_riem_prop}, there is a curve, say $c_0:I\to
M$, satisfying the above assumptions and such that
$c_0(0)=m_0,c'(0)=v_0$. Let now $c:I\to M$ be any curve satisfying
the assumptions of the Proposition. By Lemma~~\ref{2_point_lem},
there is $g\in\iso(M)$ such that $g(m_0)=c(0),g\cdot v_0=c'(0)$.
By Remark~~\ref{2_point_rem}, we can assume that $g$ preserves
orientation. By the uniqueness claim in
Lemma~~\ref{exi_uni_riem_prop}, we have $c(t)=g\cdot c_0(t)$.
\end{proof}
\end{prop}

We will now extend Theorem~~\ref{deriv_norm_thm} to the present
setting.

\begin{thm}            \label{Riem_deriv_norm_thm}
Let $M$ be an oriented two-point homogeneous space. Let
$c(t),\,t\in I,$ be a strongly regular curve in $M$.

\noindent 1. The functions $||c'(t)||,\dots,||c^{(n)}(t)||$
determine the curve up to an isometry of $M$.

\noindent 2. The functions $||c'(t)||,\dots,||c^{(n)}(t)||$ and
the number $\sgn(c)\in\{1,-1\}$ determine the curve up to an
orientation preserving isometry of $M$.
\begin{proof}
Let $f_i:I\to\R_+,1\le i \le n,$ be positive functions, let
$\si\in\{1,-1\}$. Let $m_0\in M$ be a particular point. Let
$v_0\in T_{m_0}M$ be a vector such that $||v_0||=||c'(0)||$.
Suppose that there is a curve $c_0:I\to M$ such that i) for $1\le
k \le n$ we have $||c_0^{(k)}(t)||=f_k(t)$; ii)
$\sgn(c_0(t))=\si$; iii) $c_0(0)=m_0,c_0'(0)=v_0$.

Let now $m\in M$ be any point; let $v\in T_mM$ be such that
$||v||=||c'(0)||$. By Lemma~~\ref{2_point_lem}, there is
$g\in\iso(M)$ such that $g\cdot m_0=m,g\cdot v_0=v$. By
Remark~~\ref{2_point_rem}, we can ensure that $g$ preserves
orientation. Set $c(t)=g\cdot c_0(t)$. If $\tc:I\to M$ is any
curve that has the same norms of the derivatives, has the same number
$\sgn(c)$, passes through the same point $\tc(0)$, and has the same tangent vector $\tc'(0)$,
then, by Corollary~~\ref{riem_deriv_norm_cor}, $\tc=c$. This
proves claim 2.

Applying orientation reversing isometries to strongly regular
curves $c(\cdot)$, we do not change the norms of their derivatives
but we flip $\sgn(c)$. Hence claim 1 follows from claim 2.
\end{proof}
\end{thm}

\subsection{Curves in non-orientable riemannian manifolds}   \label{non_orient_sub}
\hfill \break In this section, $M$ is a non-orientable riemannian
manifold.\footnote{The inequality $\dim M\ge 2$ is necessarily
satisfied.} The Frenet-Serret approach works with slight
modifications. We will consider only strongly regular curves
$c:I\to M$. Thus, we assume that the vectors
$c'(t),\dots,c^{(n)}(t)\in T_{c(t)}M$ are linearly independent.
The manifold does not impose any orientation on $T_{c(t)}M$; for
$t\in I$ we set $V(t)=T_{c(t)}M$ oriented in such a way that the
vectors $c'(t),\dots,c^{(n)}(t)$ form a positive basis. Note that
we may have parameter values $t_1,t_2\in I$ such that
$c(t_1)=c(t_2)$ but $V(t_1)\ne V(t_2)$. This happens when the
curve $c(\cdot)$ passes through a point, say $m\in M$, more
than once, inducing opposite orientations on the tangent space $T_mM$.

From now on, $c:I\to M$ is a strongly regular curve.
Orthonormalizing the vectors $c'(t),\dots,c^{(n)}(t)\in V(t)$, 
we obtain the orthonormal frame
$e(t)=(e_1(t),\dots,e_n(t))$.
\begin{lem}        \label{curv_nonorien_lem}
The Frenet-Serret frame $e(\cdot)$ satisfies
equation~~\eqref{riem_frenet_eq}. The curvature functions
$\ka_1(t),\dots,\ka_{n-1}(t)$ are positive.
\begin{proof}
The argument that we used to prove
equation~~\eqref{riem_frenet_eq} applies verbatim. By
construction, the vectors $e_i(t)$ satisfy
equation~~\eqref{span_eq} for $1\le i \le n$. This implies the
positivity of all the curvatures.
\end{proof}
\end{lem}
\begin{rem}        \label{curv_nonorien_rem}
{\em

In particular, $e_1(t),\dots,e_n(t)$ is a positive orthonormal
basis in $V(t)$.

}
\end{rem}

\medskip

Lemma~~\ref{curv_nonorien_lem} and Remark~~\ref{curv_nonorien_rem}
allow us to extend the material in sections~~\ref{orient_sub}
and~~\ref{2-point_sub} to non-oriented riemannian manifolds. The
following theorem states the main claims. In order to prove them,
it suffices to repeat verbatim the proofs of homologous claims in
sections~~\ref{orient_sub} and~~\ref{2-point_sub}, invoking
Lemma~~\ref{curv_nonorien_lem} and
Remark~~\ref{curv_nonorien_rem}.

\begin{thm}        \label{nonoren_main_thm}
Let $c:I\to M$ be a strongly regular curve in a non-orientable
riemannian manifold. Let $c'(t),\dots,c^{(n)}(t)$ be the tangent
vectors and let $\ka_1(t),\dots,\ka_{n-1}(t)$ be the curvatures.
Then the following claims hold.

\noindent 1. For $1 \le r \le n-1$ we have the identities
\begin{equation}    \label{non_orien_curv_eq}
\ka_r=\frac{\vol(c',\dots,c^{(r-1)})\vol(c',\dots,c^{(r+1)})}{\vol(c',\dots,c^{(r)})^2}||c'||^{-1}.
\end{equation}

\noindent 2. Let $f_1,\dots,f_n$ be positive functions on
$I\subset\R$. Assume that $0\in I$. Let $m_0\in M$ and $v_0\in
T_{m_0}M$ be such that $||v_0||=f_1(0)$. Then there exists at most
one strongly regular curve $c:I\to M$ such that i)
$c(0)=m_0,c'(0)=v_0$; ii) for $t\in I$ we have
$||c'(t)||=f_1(t),\dots,||c^{(n)}(t)||=f_n(t)$.

\noindent 3. Suppose that $M$ is a two-point homogeneous space.
Then the functions $||c'(\cdot)||,\dots,||c^{(n)}(\cdot)||$
determine the curve up to an isometry of $M$.
\end{thm}

In conclusion we note that there are non-orientable two-point
homogeneous spaces, e. g., the even-dimensional real projective
spaces.

\medskip

\noindent{\bf Acknowledgements.} While working on the project, the author
enjoyed the hospitality of mathematical establishments at the following institutions:    
Albert-Ludwigs-Universit\"at in Freiburg im Breisgau, FIM, ETH in Zurich, and the 
University of California in Los Angeles. The author has
presented some of the results in this article at the Geometry Oberseminar
in Albert-Ludwigs-Universit\"at. The participants' feedback is greatly appreciated. 
The work was partially supported by MNiSzW grant NN201384834.






\end{document}